\documentclass[twoside]{article}
\usepackage{amsthm,amstext,amscd,latexsym,amsfonts,amssymb}
\def\note#1{}
\def\H{\mathcal{H}}

\def\0n{_{(0)}}
\def\1n{_{(1)}}
\def\2n{_{(2)}}
\def\3{_{(3)}}

\begin{document}
\title{\large\bf QUANTUM-CLASSICAL INTERACTIONS\\
\large\bf AND GALOIS TYPE EXTENSIONS
\thanks{This work was partially sponsored by the
Polish Committee for Scientific Research (KBN) under Grant No.\ 5P03B05620.}}
\author{\it W\l adys\l aw Marcinek\\
\it Institute of Theoretical Physics, University
\it of Wroc{\l}aw,\\
\it Pl.\ Maxa Borna 9, 50-204  Wroc{\l}aw,\\
\it Poland}
\date{}
\maketitle

\begin{abstract}
An algebraic model for the relation between a certain classical particle system and
the quantum environment is proposed. The quantum environment is described by the
category of possible quantum states. The initial particle system is represented by an
associative algebra in the category of states. The key new observation is that
 particle interactions
with the quantum environment can be described in terms of Hopf-Galois theory. This
opens up a possibility to use quantum groups in our model of particle interactions.
\end{abstract}

\section{Introduction}

The study of highly organized structures of matter leads to the investigation of some
non-standard physical particle systems and effects. The fractional quantum Hall effect
provides an example of a system with a well-defined internal order
\cite{tsg,eha,hal,zee}. Other interesting structures appear in the so-called
$\frac{1}{2}$-electronic magnetotransport anomaly \cite{jai,jai2,jai3,dst}, high
temperature super-conductors or laser excitations of electrons. In these cases,
anomalous behaviour of electrons occurs. An example is also given by the  concept of
statistical-spin liquids (see \cite{bys} and references therein).

It seems interesting to develop an algebraic approach to the unified description of
all these new structures and effects. To this purpose,
 it is natural to assume that
the whole world is divided into two parts: a classical particle system and its quantum
environment. The classical system represents the observed reality,  particles that
really exist. The quantum environment represents all quantum possibilities that can
become
 part of  reality in the future \cite{haa}.  The goal of this paper is to
sketch a proposal of an algebraic model to describe  interactions responsible for the
appearance of the aforementioned highly organized structures. This model is based on
a general
 algebraic formalism   of Hopf algebras and Galois
extensions of rings \cite{mon}.

Our construction can be described in two steps. The first step concerns the
transformation  of the initial particles  under interaction into composite systems
consisting of quasi-particles and quanta. Such systems represent possible results of
interactions \cite{qstat,com,gqsd,gsd}. In the second step, we describe the algebra
of realizations of quantum possibilities. This step is connected with the
construction of an algebra extension and with a `decision' which possibility can be
realized and which one cannot. The problem  how such `decisions' are made was solved
in \cite{qsym} with the help
 of quantum commutativity and generalized Pauli exclusion principle.
 Our approach is based on the previously
developed concept of particle systems with generalized statistics and quantum
symmetries \cite{castat,qstat,com,gqsd,gsd,gsi}.

The paper is organized as follows. In Section 2, we propose our general model within
the framework of Hopf algebras and Galois extensions of rings. Then, in the subsequent
section, we review Hopf-algebraic generalities. This recalls mathematical concepts
employed in our proposal and allows us to specialize in the final section to the
setting already appearing in some physical models. Since there are many interesting
finite quantum groups related to spin coverings (e.g., \cite{c-a,dr}) and the
appearance of quantum symmetry in physics is more and more pronounced (e.g.,
\cite{ck}), our hope is that our general model will help us to understand some
physical phenomena that cannot be adequately described by earlier methods.

\section{The main idea}

Let us consider a system of charged particles interacting with an external quantum
environment. We assume that every charge is equipped with the ability to absorb and
emit quanta of a certain nature. A system that contains a charge and a certain
 number of
quanta as a result of interaction with the quantum environment is said to be a {\it
dressed particle} \cite{top,wmq}. A particle dressed with a single quantum is a
fictitious particle called {\it quasi-particle}. Our model is based on the assumption
that every charged particle  transforms under interaction into a composite system
consisting of quasi-particles and quanta \cite{com}. This system represents possible
results of interactions. Note that the process of absorption of quanta by a charged
particle should be described as the creation of quasi-particles, whereas the emission
as the annihilation of quasi-particles.

In our model the quantum environment is represented by a tensor category $\mathcal{C}
= \mathcal{C}(\otimes, {\bf k})$ with duals \cite{castat}. All possible physical
processes are represented as arrows of the category $\mathcal{C}$. If
$f:\mathcal{U}\rightarrow\mathcal{V}$ is an arrow from $\mathcal{U}$ to
$\mathcal{V}$, then the object $\mathcal{U}$ represents physical objects before
interactions and $\mathcal{V}$ represents possible results of interactions. Different
objects of the category $\mathcal{C}$ describe physical objects of different nature,
charged particles, quasi-particles or different species of quanta of an external
field, etc. If $\mathcal{U}$ is an object of $\mathcal{C}$ representing particles,
quasi-particles or quanta, then the object $\mathcal{U}^{\ast}$ corresponds to
anti-particles, or quasi-holes or dual fields, respectively. In the same fashion, if
$\mathcal{U}$ represents charged particles and $\mathcal{V}$ describes certain
quanta, then the product $\mathcal{U}\otimes\mathcal{V}$ encodes a composite system
containing both particles and  quanta. An arrow
$\mathcal{U}\rightarrow\mathcal{U}\otimes\mathcal{V}$ means an interaction causing
the passage from a single particle state to
 a composite quantum system.
 Thus the arrow $\mathcal{U}\rightarrow\mathcal{U}\otimes\mathcal{V}$
describes a process of absorption. Much in the same way, we conclude that the arrow
$\mathcal{U}\otimes\mathcal{V}\rightarrow\mathcal{U}$ describes a process of emission.

In our approach a unital and associative algebra $\mathcal{A}$ in the category
$\mathcal{C}$ represents the classical states of a system. The multiplication
$m:\mathcal{A}\otimes\mathcal{A}\rightarrow\mathcal{A}$ is a morphism in this
category representing the creation  of a single object of reality
 from a composite system of objects of the same species.  Quanta are
encoded in a  finitely generated coquasitriangular Hopf algebra $\H$. Quasi-particles
are described by a new algebra $\mathcal{A}^{\it ext}$, which is an extension
 of $\mathcal{A}$. Interactions are described by a right
action and  coaction of $\mathcal{H}$ on the algebra $\mathcal{A}^{\it ext}$. \note{
Obviously, the right action and  coaction of $\H$ on $\mathcal{A}^{\it ext}$ must be
compatible, so that we assume that $\mathcal{A}^{\it ext}$ is a right
$\mathcal{H}$-Hopf module \cite[p.15]{mon}. } It is  natural to assume that the
algebra $\mathcal{A}$ is invariant and coinvariant with respect to the action and
coaction of $\H$, respectively, i.e., $ \mathcal{A}=(\mathcal{A}^{{\it ext}})^{\H}
=(\mathcal{A}^{\it ext})^{\it co \H}. $ (Here $(\mathcal{A}^{\it ext})^{\H}$ is the
set of $\H$-invariants and $(\mathcal{A}^{\it ext})^{\it co \H}$ the set of
$\H$-coinvariants.) \note{ We can consider a left action and  coaction similarly, and
assume that the left and right structures are compatible. Then we obtain a Hopf
bimodule. We assume that the right (co)action can be transformed into the left one by
quantum commutativity. }

We would like to represent the interaction  of a charged particle with external
quanta as a process of creation or annihilation of quasi-particles. A composite
system of  quasi-particles and quanta is described by a tensor product
$\mathcal{A}^{\it ext}\otimes\mathcal{H}$ representing all possible quantum
configurations coming as a result of the quantum absorption process.
 On the other hand,  a
composite system of two quasi-particles (related to the same particle)
 is described by a tensor product
$\mbox{$\mathcal{A}^{\it ext}\,_{\mathcal{A}}\!\otimes\mathcal{A}^{\it ext}$}$. When
$\H$ is a group-ring Hopf algebra ${\bf k}G$, our model assumes that an element of
$G$ can be understood as a specific charge characterizing an internal degree of
freedom of a quasiparticle. We also assume that these charges are additive.
Mathematically, this means that the algebra $\mathcal{A}^{\it ext}$ is  $G$-graded,
and that this grading is strong.
 As explained in the subsequent
sections, the strongness of the $G$-grading of ${\cal A}^{ext}$
 is known to be equivalent to the bijectivity of the canonical
map
\begin{equation}
\beta :\mathcal{A}^{\it ext}\;_{\mathcal{A}}\otimes\mathcal{A}^{\it
ext}\rightarrow\mathcal{A}^{\it ext}\otimes \mathcal{H}.
\end{equation}
The bijectivity of this map means that the coaction $\mathcal{A}^{\it ext}\rightarrow
\mathcal{A}^{\it ext}\otimes\H$ is Galois. ($\mathcal{A}^{\it ext}$ is a Hopf-Galois
$\H$-extension of $\mathcal{A}$.)

An advantage of the above Galois condition is that it makes sense for an arbitrary
Hopf algebra $\H$ and does not force $\mathcal{A}^{\it ext}$ to be a crossed-product
algebra $\mathcal{A}^{\it ext} \rtimes\H$ \cite[p.101]{mon}. Thus, if we think of
$\H$ as `the group algebra of a quantum group', we have a rather general mathematical
formalism capable of describing quasi-particles with the possible charges that are
labeled by `the elements of  a quantum group' and additive according to the
multiplication of $\H$. This way the Galois condition corresponds to the additivity
of charges. \note{ In this way, $\mathcal{A}^{\it ext}$ should be at the same time a
Hopf bimodule and Galois extension. }

\note{ If the extended algebra $\mathcal{A}^{\it ext}$ is also a Hopf bimodule with a
non-trivial quantum commutative multiplication, then we say that $\mathcal{A}^{\it
ext}$ is an algebra of {\it realizations}. The opposite case, when the multiplication
is degenerated (i. e. is identically equal to 0), is connected with the generalized
Pauli exclusion principle. }

\section{Quantum commutativity and Hopf-Galois extensions}

Let   $\H$ be a Hopf algebra over a ground field ${\bf k}$, and let $m, \eta,
\triangle, \varepsilon, S$ denote its multiplication,  unit, comultiplication,
counit  and  antipode, respectively. We use the following notation for the coproduct
in $\H$. If $h \in \H$, then $\triangle (h) := \Sigma h_{\1n} \otimes h_{\2n} \in \H
\otimes \H$.
 We assume that $\H$ is a coquasitriangular Hopf algebra (CQTHA) (e.g., see
\cite[p.184]{mon}). This means that $\H$ is equipped with a convolution invertible
homomorphism
 $b \in Hom(\H \otimes \H, {\bf k})$ such that
\begin{eqnarray}
&&\label{1}\Sigma b(h_{\1n}, k_{\1n}) k_{\2n} h_{\2n}
= \Sigma h_{\1n} k_{\1n} b(h_{\2n}, k_{\2n}),\\
&&\label{2}b(h, kl) = \Sigma b(h_{\1n}, k) b(h_{\2n}, l),\\
&&\label{3}b(hk, l) = \Sigma b(h, l_{\2n}) b(k, l_{\1n}),
\end{eqnarray}
for every $h, k, l \in \H$. We call such a bilinear form $b$
 a {\it coquasitriangular structure} on $\H$.

\note{ Let us briefly recall the notion of a Hopf module (see \cite{mon} for
details). If $M$ is a right $\H$-comodule with a coaction $\delta : M \rightarrow M
\otimes \H$ of some
 Hopf algebra $\H$, then the set $M^{co \H}$ of $\H$-coinvariants is
defined by the formula
\begin{equation}
M^{co \H} := \{m \in M: \delta(m) = m \otimes 1\}.
\end{equation}
A right $\H$-Hopf module is a ${\bf k}$-linear space $M$ such that\\
(i) there is a right $\H$-module action $\lhd : M \otimes \H \rightarrow M$,\\
(ii) there is a right $\H$-comodule map $\delta : M \rightarrow M \otimes \H$, \\
(iii) $\delta$ is a right $\H$-module map, i.e.,
\begin{equation}
\Sigma \ (m \lhd h)\0n \otimes (m \lhd h)\1n = \Sigma \ m\0n \lhd h\1n \otimes m\1n \
h\2n,
\end{equation}
where $\delta (m) = \Sigma \ m\0n \otimes m\1n$. The concept of left Hopf
$\H$-modules can be introduced in a similar way. If $M$ is a right and left $\H$-Hopf
module and the left and right Hopf-module structures are compatible, then it is said
to be an $\H$-Hopf bimodule. If $\H$ is a CQTHA, then every right $\H$-Hopf module is
an $\H$-Hopf bimodule. In this case, the right and left coactions coincide (up to a
constant). }

Another ingredient of our model is the concept of quantum commutativity
\cite{cowe,qsym}. Let $\mathcal{A}$ be a unital and associative algebra and $\H$ be
a  Hopf algebra. If $\mathcal{A}$ is a right $\H$-comodule such that the
multiplication map $m : \mathcal{A}\otimes\mathcal{A}\rightarrow\mathcal{A}$ and the
unit map $\eta : {\bf k}\rightarrow\mathcal{A}$ are $\H$-comodule maps, then we say
that it is a right $\H$-comodule algebra. The algebra $\mathcal{A}$ is said to be
{\it quantum commutative} with respect to the coaction of $\H$ and its
coquasitriangular structure $b$ if an only if we have the relation
\begin{equation}\label{qc}
\begin{array}{c}
a \ b = \Sigma \ b(a_{(1)}, b_{(1)}) \ b_{(0)} \ a_{(0)} .
\end{array}
\end{equation}
Here $\rho (a) = \Sigma a_{(0)}\otimes a_{(1)}\in\mathcal{A} \otimes \H$, and $\rho
(b)=\Sigma b_{(0)}\otimes b_{(1)} \in \mathcal{A} \otimes \H$ for every $a, b \in
\mathcal{A}$. The Hopf algebra $\H$ is called a {\it quantum symmetry} of
$\mathcal{A}$.

Finally, let us recall the definition of a Hopf-Galois extension. An algebra
extension $\mathcal{A}^{\it ext}$ of $\mathcal{A}$ such that it is a right
$\mathcal{H}$-comodule algebra and $\mathcal{A}$ is its coinvariant subalgebra
\begin{equation}
\begin{array}{c}
{\mathcal A}\equiv({\mathcal A}^{\it ext})^{\it co{\mathcal H}} :=\{a \in{\mathcal
A}^{\it ext}:\delta(a)=a\otimes 1\}
\end{array}
\end{equation}
is said to be an $\mathcal{H}$-extension. If in addition the map $\beta
:\mathcal{A}^{\it ext}\;_{\mathcal{A}}\otimes\mathcal{A}^{\it
ext}\rightarrow\mathcal{A}^{\it ext}\otimes \mathcal{H}$ defined by
\begin{equation}
\begin{array}{c}
\beta (a\;_{\mathcal{A}}\otimes b) := (a \otimes 1)\delta (b) \label{gal}
\end{array}
\end{equation}
is bijective, then the $\mathcal{H}$-extension is called Hopf-Galois. If
$\mathcal{A}^{\it ext}$ is a Hopf-Galois $\mathcal{H}$-extension, then there is also a
bijection
\begin{equation}
\begin{array}{c}
\beta^n:\underbrace{{\mathcal A}^{\it ext}\;_{{\mathcal A}}\otimes\cdots\;_{{\mathcal
A}}\otimes{\mathcal A}^{\it ext}}_{n+1}\leftrightarrow{\mathcal A}^{\it ext}
\otimes\underbrace{{\mathcal H}\otimes\cdots\otimes{\mathcal H}}_{n}
\end{array}
\end{equation}
 given by
\begin{equation}
\begin{array}{c}
\beta^n := (\beta\otimes id)\circ\cdots\circ(id\;_{{\mathcal A}}\otimes\beta\otimes
id)\circ(id\;_{{\mathcal A}}\otimes\beta).
\end{array}
\end{equation}
In our physical interpretation, the one-to-one correspondence $\beta^n$ means that
the $n$-th $_{{\mathcal A}}\otimes$-tensor product representing a composite system of
$n$ quasi-particles also corresponds to
 a system of a single quasi-particle and $n$ quanta.

\section{Strongly \boldmath$G$-graded quantum-commutative algebras}

 Recall first that the group algebra ${\bf k}G$ is a
Hopf algebra for which the comultiplication, the counit, and the antipode are given
by the formulae
$$
\begin{array}{ccc} \triangle (g) := g \otimes g,&\varepsilon(g)
:= 1,&S(g) := g^{-1},
\end{array}
$$
respectively. The coquasitriangular structure on ${\bf k}G$ is given by a commutation
factor $b : G\times G\rightarrow{\bf k}\setminus\{0\}$ \cite{mon,sch,WM4,zoz}, and
the category of right $\H$-comodules is equivalent to the category  of $G$-graded
vector spaces.

Next,  assume that an algebra $\mathcal{A}^{ext}$ is an object of this category. This
means that it is a $G$-graded algebra. Now we come to the crucial theorem
\cite[p.126]{mon} stating that, for an arbitrary $G$-graded algebra and ${\bf
k}G$-coaction compatible with the grading ($\rho(a)=a\otimes g$ for $a\in
\mathcal{A}^{ext}_g$), the coaction is {\em Galois} if and only if the algebra is
{\em strongly} $G$-graded. The latter means that
\begin{equation}
\begin{array}{ccc}
{\mathcal{A}}^{\it ext}=\bigoplus_{g\in G}{\mathcal{A}}^{\it ext}_g,&
{\mathcal{A}}^{\it ext}_g{\mathcal{A}}^{\it ext}_h = {\mathcal{A}}^{\it ext}_{gh},&
{\mathcal{A}}^{\it ext}_e\equiv{\mathcal{A}},
\end{array}
\end{equation}
where $e$ is the neutral element of $G$.

As an example, let us  consider a $G$-graded $b$-commutative $\Bbb C$-algebra
$\mathcal{A}^{ext}$ with the so-called standard gradation \cite{WM4}. This means that
we take as the {\em strongly} grading group $Z^N := Z\oplus...\oplus Z$ and assume
\begin{equation}\label{standard}
b (\xi^i , \xi^j)=:b^{ij} = (-1)^{\Sigma_{ij}} q^{\Omega_{ij}}. \label{comd}
\end{equation}
Here $\xi^i := (0, \ldots, 1, \ldots, 0)$ ($1$  on the $i$-th place) is the set of
generators of $Z^N$, $\Sigma := (\Sigma_{ij})$ and $\Omega := (\Omega_{ij})$ are
integer-valued matrices such that $\Sigma_{ij} = \Sigma_{ji}$ and $\Omega_{ij} = -
\Omega_{ji}$, and $q \in {\Bbb C} \setminus \{0\}$ is a parameter \cite{zoz}. Since
our Hopf algebra is a group ring, the equation (\ref{1}) for the bilinear form $b$ is
automatically satisfied, whereas the equations (\ref{2})-(\ref{3}) uniquely determine
$b$ once we set its value on the generators. The convolution-invertibility of $b$
follows from the fact that $b^{ij}$s are always non-zero. (Notice that
$b(\xi^i,(\xi^j)^n)=(b^{ij})^n$, so that for  $q = exp(\frac{2 \pi i}{n})$ and $n$
even, the grading group $Z^N$ can be reduced to $Z_n \oplus...\oplus Z_n$.) Combining
(\ref{qc}) with (\ref{standard}), we obtain the following quantum commutativity
relations:
\begin{equation}
a_{\xi^i}\;a_{\xi^j}=b^{ij}\;a_{\xi^j}\; a_{\xi^i},\;\;\;\mbox{where}\;\;\;
a_{\xi^i}\in\mathcal{A}^{ext}_{\xi^i},\;\;\; a_{\xi^j}\in\mathcal{A}^{ext}_{\xi^j}.
\end{equation}
 It is the behaviour of $b^{ij}$ that determines whether we obtain
a system with the $q$-statistics, or  Fermi statistics and the Pauli exclusion
principle, or whether we obtain bosons. On the other hand, the {\em strong} gradation
ensures that the internal degrees of freedom of a quasi-particle are labeled by
charges ($N$-tuples of integers), and that these charges are {\em additive}.

{\bf Acknowledgments.} It is a pleasure to thank Cezary Juszczak for his help with
typesetting this article.

\end{document}